\documentclass[12pt]{article}
\usepackage{epsf}

\usepackage{amsmath, amsfonts,amssymb,amsthm,graphicx,color} 
\newtheorem{theorem}{Theorem}[section]
\newtheorem{proposition}[theorem]{Proposition}
\newtheorem{definition}[theorem]{Definition}
\newtheorem{remark}[theorem]{Remark}
\newtheorem{assumption}[theorem]{Assumption}
\newtheorem{example}[theorem]{Example}

\newtheorem{corollary}[theorem]{Corollary}
\newtheorem{lemma}[theorem]{Lemma}
\newtheorem{notation}[theorem]{Notation}
\newtheorem{acknowledgements}[theorem]{Acknowledgements}

\def\A{\= A}
\def\a{\= a}

\def\u{\= u}

\def\ndi{\={\i}}

\def\suls{{\'Sulvas\u tras}}

\date{}

\def\build#1_#2^#3{\mathrel{\mathop{\kern 0pt#1}\limits_{#2}^{#3}}}

\def\smallsquare{\vbox{\hrule\hbox{\vrule height 1 ex\kern 1
ex\vrule}\hrule}}

\begin{document}
      
\begin{center}
{\Large{\bf Mensuration of Quadrilaterals in the L\ndi l\a vat\ndi}}

\medskip
S.G. Dani
\end{center}

\begin{abstract}
Mensuration with quadrilaterals had received attention in the Siddh\a nta 
tradition at least
since Brahmagupta. However, in Bh\a skarac\a rya's {\it L\ndi l\a vat\ndi} we 
come across some
distinctively new features. In this paper an attempt will be made to put the 
development in historical perspective. 
\end{abstract}

\medskip
A systematic study of the topic of mensuration of quadrilaterals 
 in Indian mathematics goes back at least to the 
{\it Br\a hmasphu\d tasiddh\a nta} (628 CE) of 
Brahmagupta (born in 598 CE); in certain special cases, such 
as isosceles trapezia, some familiarity is found in the 
{\it \suls} from around the middle of the first millennium BCE 
(see Sen and Bag 1983\footnote{Sen and Bag 1983 : S.N. Sen and A.K. Bag, The {\it \suls}, 
Indian National 
Science Academy, New Delhi, India, 1983}).
With regard to mensuration of quadrilaterals Brahmagupta is
well-known for the formula in the s\u tra:

{\small

\begin{center}
{\tt sth\u laphala\d m tricaturbhujab\a hupratib\a huyogadalagh\a ta\d h |

\smallskip
bhujayog\a rdhacatu\d stayabhujonagh\a t\a t pada\d m s\u k\d sma\d m ||}
\end{center}
}
\smallskip 

~\hfill {\it  Br\a hmasphu\d tasiddh\a nta} (XII - 21)

\medskip
Traditionally  the s\u tra has been understood, by ancient mathematicians
following Brahmagupta (I shall dwell more on this later), 
as well as modern commentators broadly as follows: 

\begin{quote}
{\it The gross area of a triangle or quadrilateral is the product of half the 
sum of the opposite sides. The exact area is the square-root of the product 
of the four sets of half the sum of the sides (respectively) 
diminished by the  sides.} 

\end{quote}

This translation is taken from Saraswati 1979\footnote{Saraswati 1979 : T.A. 
Saraswati Amma, {\it Geometry in Ancient and Medieval India}, Motilal Banarasidas 
Publishers, Delhi, India,  1979, reprinted in 1999.}(p. 88), and the
translations in Colebrooke  
1817\footnote {Colebrooke 1817 : H.T. Colebrooke, {\it Algebra, arithmetic and 
mensuration, from the original 
Sanscrit, of Brahmegupta and Bhaskara}, Translation, John Murray, Albemarle 
Street, London, UK, 1817.} and Plofker 2009\footnote{Plofker 2009 : Kim Plofker, 
{\it Mathematics in India}, Princeton University Press, Princeton, NJ,
USA, 2009. xiv+357 pp.} 
 also correspond to it. In all of these, 
in particular, {\it tricaturbhuja} is interpreted as referring to 
the formula being applicable to triangles and quadrilaterals (independently). 
For the case of the triangle this is the well-known formula known after 
Heron of Alexandria (1st century~CE), and in this context the quadrilateral 
version is 
referred to as ``Brahmagupta's generalization'' (see, for instance, 
Plofker 2009, p. 144). 
The general version, which in modern notation may be stated as
 $$ \ \ \ \ \ \ \ \ \ \ \ A=\sqrt {(s-a)(s-b)(s-c)(s-d)}, \ \ \ \ \ \ \ \ \ \ \ (*)$$ where 
$a,b,c,d$ are the sides of the quadrilateral, $s$ is half
the perimeter, and $A$ is the area,  is however correct only 
when the quadrilateral is cyclic, viz. when all 
the four vertices lie on a common circle; this condition holds for 
quadrilaterals like rectangles and isosceles trapezia,  but not in general, 
e.g. for rhombuses 
with unequal diagonals. The general perception in the context of the
interpretation has been that 
the author 
somehow omitted to mention the condition, though actually aware of 
it, with the latter being confirmed by the fact that he is noted to 
apply  it only for cyclic quadrilaterals. 

It has however recently been argued in 
Kichenassamy 2010\footnote{Kichenassamy 2010 :
S.~Kichenassamy, Brahmagupta's derivation of the area of a 
cyclic quadrilateral, Historia Math. 37 (2010), 28 - 61.}
that actually  the term {\it tricaturbhuja} was in fact used by Brahmagupta to 
mean 
a cyclic quadrilateral (and {\it not} triangle and/or quadrilateral). Thus 
Brahmagupta means to state 

\begin{quote}
{\it The gross area of a {\bf cyclic quadrilateral} is the product of half the 
sum of the opposite sides, and  the square-root of the product 
of the four sets of half the sum of the sides (respectively) 
diminished by the  sides is the exact area. } 

\end{quote}

A bit of a hint towards that the traditional  
interpretation may not be right is contained in 
the fact that at the only other place where the term {\it tricaturbhuja} 
occurs in {\it Br\a hmasphu\d tasiddh\a nta}, in s\u tra (XII - 27),
(and it is not known to occur in earlier or later ancient texts) the result 
involved (relating to the circumradius) is stated first for triangles, 
separately, and then for 
 ``{\it tricaturbhuja}''s which indicates that the latter should in fact have 
four sides (see Kichenassamy 2013\footnote{Kichenassamy 2013 : S. 
Kichenassamy, Textual 
analysis of ancient Indian mathematics, 
Ga\d nita Bh\a rat\ndi\ 33 (2011), 15 - 28 (2013).}).
 The arguments in Kichenassamy 2010,
go well beyond that, with the author providing a detailed discussion on 
the issue, 
including on how Brahmagupta would have arrived at the formula, and how it
incorporates in a natural way the hypothesis that the quadrilateral is 
cyclic. According to Kichenassamy,  Brahmagupta while pursuing his study of 
triangles dealt with the circumcircle, described in particular 
a formula for the circumradius, and along the line of thought considered 
quadrilaterals formed by the triangle and a point on the circumcircle, which
motivated 
the term {\it tricaturbhuja}.\footnote{In a recent paper (not yet published) 
P.P. Divakaran proposes a somewhat different scenario for the 
development of ideas in Brahmagupta's work and the genesis of the term 
{\it tricaturbhuja}, which nevertheless discounts the traditional 
interpretation of the term mentioned earlier. }  

Unfortunately, the theory developed by Brahmagupta did not go down the line 
of later mathematicians in India with proper understanding. It may be 
worthwhile to recall the following in this respect. Let us consider the 
works of 
\'Sr\ndi dhara, the author of  {\it P\a t\ndi ga\d nita}\footnote 
{\'Sr\ndi dhara, {\it P\a t\ndi ga\d nita} : {The P\a t\ndi ga\d nita of
 \'Sr\ndi dharac\a rya 
with an ancieṇt Sanskrit commentary}, ed \& tr, Kripa Shankar 
Shukla, Hindu 
Astronomical and Mathematical Texts Series No.~2,
Department of Mathematics and Astronomy, Lucknow University, Lucknow, India, 1959.} 
and  {\it Tri\'satik\a}, and  Mah\a v\ndi ra who authored  {\it Ga\d nitas\a 
rasa\. ngraha}\footnote{Mah\a v\ndi ra,{\it Ga\d nitas\a rasa\. ngraha} : 
{Mah\a v\ndi ra, \it Ga\d nitas\a rasa\. ngraha}, ed \&
tr Dr. (Mrs.) Padmavathamma, Publ. by Sri Hombuja Jain Math,
Hombuja, Shimoga District, Karnataka, India, 2000.},  two prominent 
authors\footnote{As noted by K.S. Shukla in the introduction to his edition of 
\'Sr\ndi dhara's {\it P\a t\i ga\d nita},
\'Sr\ndi dhara's works are cited by many later authors. On the other 
hand Mah\a v\ndi ra's  {\it Ga\d nitas\a 
rasa\. ngraha} apparently enjoyed 
the status of a textbook in many parts of South India for nearly 
three centuries, until the arrival of Bh\a skar\a c\a rya's 
{\it L\ndi l\a vat\ndi}, as noted by Balachandra Rao, in his review  
(in Ga\d nita Bh\a rat\ndi Vol. 35 (2013), p.~167) of the book 
{\it \'Sr\ndi\ R\a j\a ditya' Vyavah\a raga\d nita} edited and translated 
by Padmavathamma, Krishnaveni and K.G. Prakash.} from the intervening 
period between Brahmagupta 
and Bh\a skara.  Mah\a v\ndi ra is known 
to be from around 850~CE. Concerning \'Sr\ndi dhara there has been a 
controversy 
among scholars over his period, and in particular over whether he was 
anterior or posterior to Mah\a v\ndi ra, but it now seems to be agreed 
that he
is from the 8th century.\footnote{In \S 5 of  Shukla's introduction 
to his edition of \'Sr\ndi dhara's {\it P\a t\ndi ga\d nita}, one finds a 
detailed discussion on this issue, 
concluding with his own verdict that \'Sr\ndi dhara  
``lived sometime between Mah\a v\ndi ra (850) and \A ryabha\d ta II (950)''. 
Saraswati 1979  expresses skepticism in this 
respect (See page 10; her wording is ``\'Sr\ndi dhara is probably 
earlier than 
Mah\a v\ndi ra though K.S. Shukla places him betewen 850 and 950 A.D.'').
S.D. Pathak, in his paper {\it \'Sr\ndi dhara's time and works}, Ga\d nita 
Bh\a rat\ndi, Vol 25 (2003), 146-149, which was published posthumously 
but was actually written before Shukla's edition of {\it P\a t\ndi ga\d nita} 
was published, 
argues in favour of \'Sr\ndi dhara being earlier, and, notably, in a special 
note following the article the Editor R.C. Gupta, who had himself also 
discussed the issue in an earlier paper,  mentions ``But now 
Dr. Shukla himself accepts the priority of \'Sr\ndi dhara over Mah\a v\ndi ra
(personal discussions)''. Also,  Shefali Jain in her  recent thesis 
{\it \A carya \'Sr\ndi dhar evam ac\a rya 
Mah\a v\ndi ra 
ke Ganiteeya avadano ka tulanatmak adhyayan} (Ph.D. Thesis, Shobhit 
Vishwavidyala Meerut, India, 2013.)  mentions that at the end of a manuscript of 
{\it Ga\d nitas\a rasangraha} in the Royal Asiatic Society, Bombay (MS No. 
230)  one finds the statement ``{\it kram\a dityukta\d m 
\'Sr\ndi dhar\a carye\d na bhadra\d m bh\u y\a t} '' which shows that
\'Sr\ndi dhara preceded Mah\a v\ndi ra. In the light of earlier literature
 she assigns 750 CE as the year around when he would have flourished.} 

In {\it P\a t\ndi ga\d nita}, on the issue of areas of quadrilaterals 
\'Sr\ndi dhara first gives a formula for the areas of trapezia (the usual 
one), in s\u tra 115 which is then complemented, in s\u tra 117 (see 
\'Sr\ndi dhara, {\it P\a t\ndi ga\d nita}, ed. K.S. Shukla, p. 175), with the 
following:

{\small
\begin{center}
\tt bhujayutidala\d m caturdh\a\ bhujah\ndi na\d m tadvadh\a tpada\d m ga\d nita\d m |

\smallskip
sad\d r\'s\a samalamb\a n\a masad\d r\'salambe vi\d samabahau ||
\end{center}

\ \hfill {\it P\a t\ndi ga\d nita} (117)
 }

\begin{quote}
Set down half the sum of the (four) sides (of the quadrilateral) in four 
places, (then) diminish them (respectively) by the (four) sides (of the 
quadralateral), (then) multiply (the resulting numbers) and take the square
root (of the product): this gives the area of the quadrilaterals 
having (two or more) equal sides but unequal altitudes and also of 
quadrilaterals having unequal sides and unequal altitudes. - (Translation 
from \'Sr\ndi dhara, {\it P\a t\ndi ga\d nita}, ed \& tr K.S. Shukla.) 
 
\end{quote}

Thus \'Sr\ndi dhara is seen to give formula ($*$) for the area of 
any quadrilateral, 
dwelling elaborately on the generality, without realising that it is not
true in that generality.
Interestingly, unlike in the case of the formula for trapezia, no examples 
are discussed 
to illustrate the general formula.  This suggests in a way that the s\u tra is 
included in the spirit of recording and passing on a piece of traditional 
knowledge in which the author espouses no direct interest; this is a kind
of situation in which mathematicians are prone to let down their guard!
His treatment 
in {\it Tri\'satik\a} is also along the same lines (See Saraswati 1979, p.~92). 

Mah\a v\ndi ra states  the result under discussion as follows: 

{\small
\begin{center}
{\tt bhujayutyardhacatu\d sk\a dbhujah\ndi naddh\a titatpadam s\u k\d sma\d m |}
\end{center}
}

\ \hfill {first half of \it Ga\d nitas\a rasa\. ngraha} (VII - 50)
 
\begin{quote}
Four quantities represented (respectively) by half the sum of the sides as 
diminished by (each of) the sides (taken in order) are multiplied together 
and the square root (of the product so obtained) gives the minutely accurate 
measure (of the area of the figure). - (Translation from Mah\a v\ndi ra, 
{\it Ga\d nitas\a rasa\. ngraha}, ed Padmavathamma).

\end{quote}

Here again formula ($*$) is stated unconditionally for any quadrilateral.  
The second half of the above mentioned verse is the usual formula for the 
area of a trapezium, as the product of the perpendicular height with half
the sum of the base and the opposite side, mentioning also a caviat that 
it does not hold for 
a {\it vi\d samacaturasra}. In verses 51--53 following the s\u tra as above, 
Mah\a v\ndi ra asks to compute the areas of tringles with given lengths 
for their sides, presumably meant to be done using the first part of 
verse~50. verse~54 describes a formula for the diagonal of a 
quadrilateral\footnote 
{This formula also goes back to Brahmagupta and is valid only for cyclic 
quadrilaterals, but is stated here unconditionally.} and in verses 55 to 57 
the author asks to compute diagonals and areas of quadrilaterals which are
isosceles trapezia; since reference to diagonals is also invoked it is not 
clear whether the computation is meant to be done using the general form of  
verse~50, namely formula ($*$), or by first computing the diagonal.  
As a whole the treatment suggests a lack of interest (perhaps coupled with disbelief)
in the general case of the formula. 

Similarly, \'Sr\ndi pati 
(11th century) also gives, in {\it Siddh\a nta\'sekhara}, formula~($*$)
unconditionally (see Saraswati 1979, p. 94). 
On the whole the 
practice of treating the expression as the formula for the area of any 
quadrilateral was so prevalent that one finds it presented as such  even in 
the 14th century work {\it Ga\d nitas\a rakaumud\ndi} of \d Thakkura 
Pher\u\ (see SaKHYa, 2009\footnote{SaKHYa, 
{\it Ga\d nitas\a rakaumud\ndi; the Moonlight of the Essence of
Mathematics, by \d Thakkura Pher\u}; edited with Introduction,
Translation, and Mathematical Commentary, Manohar Publishers, New
Delhi, India,  2009.}, p. 142). Notwithstanding the overall 
continuity of the Indian mathematical tradition, topics that were not directly 
involved in practice, in astronomy or other spheres in which mathematics was 
applied at the time, suffered neglect, and sometimes were  
carried forward without a proper understanding of what was involved. 

By the time of Bh\a skara any connection of the formula ($*$) with cyclicity 
of the quadrilateral was completely lost. Even  
awareness of cyclic quadrilaterals seems to have gone missing over a period,  
until it was resurrected in the work of N\a r\a ya\d na Pa\d n\d dita, in 
the 14th century (See Saraswati, 1979, pages 96-106, for details).    

In this overall context, as it  prevailed around the turn of the 
millennium,  \A ryabha\d ta II, who is believed to have 
lived sometime between 950 and 1100 CE (see Plofker 2009, p.~322), 
rejected ($*$) 
as the formula for  the area of a quadrilateral, ridiculing one 
who wants to find the area of a quadrilateral without knowing the 
length of a diagonal as a fool or devil ({\it m\u rkha\d h pi\'s\a co v\a}) 
(see Saraswati 1979, p. 87). This was the situation when 
Bh\a skara appeared on the scene. Though apparently guided by 
 \A ryabha\d ta~II in his treatment in respect of Brahmagupta's theorem
(see Saraswati, 1979, p.~94) Bh\a skara took an entirely different approach 
to the issue, bringing considerable clarity on the topic (even 
though he did not get to cyclic quadrilaterals). 

A closer look at the relevant portion of the {\it L\ndi l\a vat\ndi} shows 
an intense concern on the part of Bh\a skara at what he observed as 
a flaw in the ``traditional'' formula. This does not seem to have been 
adequately appreciated in the literature on the topic. One of the 
reasons for this seems to be that the standard translation cum commentaries
(E.g. Colebrooke 1993\footnote{ Colebrooke 1993 : H.T. Colebrooke, 
Translation of the {\it Lilavati}, with notes by Haran 
Chandra Banerji, Second Edition, Asian Educational Services, New
Delhi, Madras, India, 
1993.}, Phadake 1971\footnote{Phadake 1971 : N.H. Phadake, 
{\it Shrimad Bhaskaracharya k\d rta 
L\ndi l\a vat\ndi\ Punardarshan}, (in {Marathi}), Sarita Prakashan, 
Pune, India, 1971, last reprinted 2014.}, Patwardhan-Naimpally-Singh 
2001\footnote{Patwardhan-Naimpally-Singh 2001 : K.S. Patwardhan, S.A. 
Naimpally and S.L. Singh, {\it L\ndi l\a vat\ndi\ 
of Bh\a skar\a c\a rya, A Treatise of Mathematics of Vedic Tradition},
Motilal  Banarasidas Publishers, Delhi, India, 2001.}, and  Jha
2008\footnote{Jha 2008 : Pandit  
Shrilashanlal Jha, {\it Shrimad Bhaskaravirachita Lilavati} 
(in {Hindi}), Chaukhamba Prakashan, Varanasi, India, 2008.}), which have been the
chief sources for dissemination of the topic, have translated and commented
upon the s\u tras involved only individually, in a rather disjointed way, 
as a result of which a common strand that  Bh\a skara  followed in respect of 
the above seems to have been missed. Secondly, many of the commentaries, 
except Colebrooke 1993 from the above, while 
including the text of the {\it L\ndi l\a vat\ndi} do not include Bh\a skara's 
{\it V\a san\a bh\a \d sya} (explanatory annotation, in prose form, on the 
original verse text) along with it; also, 
even as they are
seen to avail of various points made in 
{\it V\a san\a bh\a \d sya}, no reference is made to the latter, which 
diffuses the overall picture even further. Colebrooke, 1993, 
 does include {\it V\a san\a bh\a \d sya} and 
also a meticulous translation of it for the most part, though 
as I shall point out below a crucial line relevant to 
the theme under discussion is missing from the
translation\footnote{Similarly some other short segments aimed at 
introducing   the subsequent verse have been omitted from the translation
  elsewhere in the text.}. 

I shall now present the part of the {\it L\ndi l\a vat\ndi} together with 
 the {\it V\a san\a bh\a \d sya} on the 
issue as above and bring out the strand of 
Bh\a skara's thinking, and concern, over the perceived flaw. 

Let me begin with Bh\a skar\a c\a rya's statement on the formula
(for reference to the verses I shall indicate  their  numbers  in 
Colebrooke~1993 and also Apte~1937\footnote{ Apte~1937 : Vinayak Ganesh Apte, 
{\it Buddhivil\a sin\ndi\ L\ndi l\a vat\ndi 
vivara\d n\a khyat\ndi k\a dvayopeta \'Shr\ndi mad Bh\a skar\a c\a ryaviracita
L\ndi l\a vat\ndi} (in {Sanskrit}), \A nand\a shram Mudra\d n\a laya, 
Pune, India, 1937.}, the standard Sanskrit sourcebook on the topic; 
the latter number will be put between square bracket, following the former,
separated by a colon. 
  The reader is cautioned that the numbers vary somewhat 
according to the reference concerned):

{\small
\begin{center}
{\tt sarvadoryutidala\d m catu\d hsthita\d m b\a hubhirvirahita\d m ca 
tadvadh\a t |

\smallskip
m\u lamasphu\d taphala\d m caturbhuje spa\d s\d tamevamudita\d m trib\a huke || 
}

\ \hfill {\it L\ndi l\a vat\ndi} - (167 : [169])

\end{center} 
}

\begin{quote}
Half the sum of all the 
sides is set down in four places; and the sides are severally subtracted. 
The remainders being multiplied together, the square root of the product 
is the area, inexact in the quadrilateral, but pronounced exact in the 
triangle. -- (The translation is taken from Colebrooke, 1993; I have added a 
comma, after ``area'', which seems to be needed for easy comprehension.) 

\end{quote}

Thus,  unlike \A ryabha\d ta II, Bh\a skara does not 
reject outright the 
formula for the quadrilateral, accepting it only for triangles. He mentions
it as exact ({\it spa\d s\d ta}) for triangles while for quadrilaterals 
he calls is ``inexact'' ({\it asphu\d ta}). This is however only the beginning. 
Detailed comments on it are to follow. 

To begin with he asks, in verse (168 : [170]), for the area to be computed 
``as told by the ancients'' ({\it tatkathita\d m yad\a dyai\d h}) for a 
quadrilateral 
with  base ({\it bh\u mi})~14, face ({\it mukha})~9, sides 13 and 
12, and perpendicular 12; it can be seen that the quadrilateral is a 
non-isosceles trapezium, formed by attaching to a rectangle with sides 9 and 12, a right 
angled triangle with sides 5, 12 and 13, along the side with length 12. 
In the {\it V\a san\a bh\a \d sya} Bh\a skara proceeds to note 
that the area given by the formula is $\sqrt {19800}$, which is ``a 
little less than 141'', while the actual area (which can be computed for 
a trapezium more directly) is 138, thus pointing 
to a contradiction.   

Though this would have sufficed for the contention that the formula is 
not correct, or accurate, Bh\a skara does not leave it at that. He embarks on 
a more detailed discussion on the issue. At the beginning of 
verse (169 : [171]), 
in the {\it V\a san\a bh\a \d sya} the \A c\a rya says:

{\small
\begin{center}
{\tt atha sth\u latvanir\u pa\d n\a rtha\d m s\u tra\d m s\a rdhav\d rtta\d m |}
\end{center}
}

\begin{quote}
Now a {\it s\u tra} of a stanza and half for looking into the grossness. (My
translation). 

\end{quote}

\noindent This line, which is significant from our point of view, is not
found in the translation in Colebrooke, 1993 (as noted earlier many other 
sources do not include the text of the {\it V\a san\a bh\a \d sya}
either). We note that ``{\it atha}'' marks commencement (of a story, 
chapter, argument etc., typically in a ceremonial way) and that 
{\it nir\u pa\d na} means ``looking into, analysis or investigation''. 
{\it sth\u latva} stands for ``grossness'' or ``coarseness''; thus Bh\a skara
is announcing here that he is taking up an analysis of the 
grossness\footnote{Like the English word adopted here for 
{\it sth\u latva}, the 
latter also has shades of meaning of unflattering variety, and one may 
wonder whether the choice of the word {\it sth\u latva} here, as against
say the noun form of the adjective {\it asphu\d ta} that was adopted 
in the original s\u tra), is deliberate.} 
 (of 
the formula). This is followed by the following
argument: 
{\small
\begin{center}
{\tt caturbhujasy\a niyatau hi kar\d nau katha\d m tatosminniyata\d m 
phala\d m sy\a t | 

\smallskip
pras\a dhitau tacchrava\d nau yad\a dyai\d h svakalpitau t\a vitaratra na 
sta\d h || 

\medskip
te\d sveva b\a hu\d svaparau ca kar\d n\a vanekadh\a\ k\'setraphala\d m 
tata\'sca  ||
}
\end{center}
}
\ \hfill {\it L\ndi l\a vat\ndi} (169-170 : [171])

\begin{quote}
The diagonals of a quadrilateral are indeterminate ({\it aniyatau});
then how could the area [confined] within them be determinate? The 
(values for) diagonals 
set down ({\it pras\a dhitau}) by the 
ancients
would not 
be the same elsewhere. For the same (choices of the) sides the diagonals have 
many possibilities
and the area would vary accordingly. (My translation).  

\end{quote}

In the {\it V\a san\a bh\a \d sya} this is further elaborated, noting that
if in a quadrilateral two opposite vertices are moved towards each other then 
the diagonal between them contracts, while the other two vertices move away
from each other and the diagonal between them elongates, and thus with sides 
of the same length there are other possible values for the diagonals.  

The next verse  in the {\it L\ndi l\a vat\ndi} raises some rhetorical 
questions: 

{\small
\begin{center}
{\tt lambayo\d h kar\d nayorvaika\d m anirdi\'syapar\d h katha\d m |

\smallskip
p\d rcchatyanitvepi niyata\d m c\a pi tatphala\d m ||}
\end{center}
}
\ \hfill  {\it L\ndi l\a vat\ndi} (171 : [173])

\begin{quote}
When none of the perpendiculars nor either of the diagonals are specified, 
how will the other values get determined? It is like asking for 
definite area, when in fact it is indeterminate. (My translation).

 \end{quote}

Then come some devastating blows:

{\small
\begin{center}
{\tt sa p\d rcchak\a h pi\'s\a co v\a\ vakt\a\ v\a\ nitar\a \d m tata\d h |

\smallskip
yo na vetti caturb\a huk\'setrasy\a niyat\a m sthiti\d m ||}

\end{center}
}
\ \hfill  {\it L\ndi l\a vat\ndi} (172 : [172])

\begin{quote} 
Such a questioner is a blundering devil ({\it pi\'s\a ca}) and worse is 
one who answers it. They do not realise the indeterminate nature of the 
area of a quadrilateral. (My translation).

\end{quote}

The scorn being deployed is reminiscent of Aryabhata~II, but here we find it
accentuated, and its scope extended to those answering the question!

Having vented his ire over the ignoramuses the \A c\a rya next sets out 
to establish the point through more concrete illustrations. Before continuing 
with it, it may be worthwhile to note the following. While the argument given 
in (169-170) is of considerable heuristic value, it is {\it not} conclusive 
from a logical point of view, since so far it has not been shown that when 
the sides are the same and diagonals vary the areas could actually be 
different. To make the argument foolproof, and to fully convince 
skeptics,  one 
needs concrete examples, 
with same four sides and different pairs of diagonals for which the areas 
{\it are verifiably different}. While he may or may not have have specifically followed 
such a train of thought,  
that is what Bh\a skara sets out to do in the following verses. 

For the illustrations he needs situations where after making alterations  
as proposed in the  {\it V\a san\a bh\a \d sya} following s\u tras 
(169-170 : [171])
it would be possible to readily compute the area (and show that it is 
different). For this purpose he considers the class of equilateral 
quadrilaterals (in which all sides are equal, also called rhombuses). He notes
a formula for the second diagonal, given the common value of the side 
and one of the diagonals: in modern notation, if $a$ is the side and $d_1$ and
$d_2$ are the diagonals then $$d_2=\sqrt {4a^2-d_1^2}.$$ He recalls also the
formula for the area, as equal to $\frac 12 d_1d_2$; the s\u tras 
(174-175) involving this include also statement of areas of rectangles 
and trapezia, but that is purely to put the situation of equilaterals in 
context - after all, the knowledge of these is already implicit in
particular in the 
problem posed in verse (169 : [170]). 

At this juncture one is in a position to illustrate the point that 
was made following (169-170 : [171]), since if we start with an equilateral 
quadrilateral with side  $a$ and a  diagonal of size $d_1$ and move the 
vertices on the two sides of that diagonal 
closer along the diagonal, then the area is 
$\frac 12 d_1 \sqrt {4a^2-d_1^2}$, with $d_1$ as a variable quantity; one can
readily see that choosing different values of $d_1$ we get rhombuses with 
the same side-lengths but different areas.
However, 
not content with this, he seeks more concrete choices for the side 
and diagonal, 
specifically with integer values. For this purpose he brings in his
knowledge of what are now called Pythagorean triples. A Pythagorean triple is a
triple of natural numbers $(l,m,n)$ such that $l^2+m^2=n^2$; by the Pythagoras 
theorem (or rather its converse, which can be deduced from the theorem itself, 
via elementary geometry) for a triangle with sides $l,m$ and $n$ where $(l,m,n)$ is a 
Pythagorean triple, the angle opposite to the side $n$ is a right angle. 
Putting 4 such right angled triangles together, along their equal sides 
adjacent to the right angles,  we get an equilateral 
quadrilateral with all sides $n$ and diagonals $2l$ and  $2m$ respectively, 
and their areas are $2lm$. 
Bh\a skara now chooses the triples $(15,20,25)$ and $(7,24,25)$ which
are seen to be Pythagorean\footnote{It may be recalled here as an aside that
Pythagorean triples have been known
in India since the time of \suls;  for 
a discussion on such triples occurring in \suls\ the reader is referred 
to: S.G. Dani, 
Pythagorean triples in \suls, Current Science 
85 (2003), 219 - 224.}, and have common value 
for the length of the hypotenuse.  Thus the construction as above yields two 
equilateral
quadrilaterals with areas $2(15\times 20)=600$ and $2(7\times 24)= 336$, 
respectively. Bh\a skara also points out that we may also consider for
comparison the square with all sides 25, in which case the area is 625, 
a yet another value for the area. 

Thus Bh\a skara adopts various means, argumentation, pressurising through 
rhetoric, as well as persuasion, to put it across to his readers that 
Brahmagupta's formula is not valid exactly for a general quadrilateral. 

Following the group of versers discussed above, there are two more problems 
concerning the area of a quadrilateral. In verse (177 : [175]) we have an 
example of (what turns out to be) a non-isosceles trapezium, for which 
again it is pointed out, in  {\it V\a san\a bh\a \d sya}, that the
area computed using ($*$) does not give the true value. 

The next verse, (178 : [176]), asks to find the area, and also the diagonals
and perpendiculars, of a quadrilateral whose sides are given as, face 51, 
base 75, left side 68 and right side 40. To a discerning reader it should 
seem puzzling that the 
\A c\a rya should ask such a question, giving only the sizes of the four 
sides, after all the painstaking endeavour to 
get it across that sizes of four sides do not determine a quadrilateral,
and in particular the area is indeterminate. It is hard to reconcile this
especially with verse  (172 : [172]), according to which one asking such 
a question is a ``{\it pi\'s\a ca}''. 
 
The spirit of what follows  however seems 
to be to explain how one should proceed in response, when such a problem 
is posed (e.g. as a challenge); since at one time the focus was 
on cyclic
quadrilaterals which were determined once the sides were given (together with 
their order) it may have been a general practice to pose questions about
quadrilaterals purely in terms of their sides (in specific order). 
In the next verse Bh\a skara recalls that if we know the 
perpendicular that determines
the (corresponding) diagonal, and knowing a diagonal determines the 
(corresponding) perpendicular, and the area; 
this is consistent with the earlier contention about the need for an 
additional assumption being necessary.    
In his treatment of the problem in the {\it V\a san\a bh\a \d sya} he then 
says ``to determine the perpendicular  we {\it assume} the diagonal
joining the tip of the left side to the base of the right side to be 
(of length) 77'' (emphasis added). 
With this choice for the diagonal the perpendicular and 
then the area of the quadrilateral are computed (adding the areas of the 
two triangles on the two sides of the diagonal as above)\footnote{Having 
assumed the diagonal to have length 77, the areas of the two 
triangles formed could be computed using ($*$), which 
Bh\a skara does consider to be exact for triangles, but the method given 
goes through
the computation of the perpendicular, and no reference is made to this other 
possibility.}; it turns out to be 3234. 
Interestingly, this is the value that one 
would get from~($*$) with the values 51, 68, 75 and 40 for the four sides!
(One may wonder whether Bh\a skara intended this, but there is no way to know.)
The reason for this agreement, from a modern perspective, is of course 
that for the above choice of the diagonal the quadrilateral is cyclic, 
and Brahmagupta's formula does apply. 

It would seem curious that the choice made 
was such that the quadrilateral is cyclic, especially when there is no 
reference to such a concept in the text. 
Also, though the fact of having to make a 
choice has been clarified, one would wonder how the choice of 77 as the 
length of the (particular) diagonal came about, especially in the context 
of its turning out to be one for which the quadrilateral is cyclic; 
it may be noted that if one starts with  
an ad  hoc choice the computations of the perpendicular
and the area involve rather complicated surds, making it unsuitable for an 
illustrative example.  All commentators have 
repeated the part about assuming the (particular) diagonal to be 77 
(generally without
reference to the  {\it V\a san\a bh\a \d sya}), but throw no 
light on the issue of what motivates the specific value that is chosen. 
It seems that 
this quadrilateral was familiar to Bh\a skara, together with the value 
for the diagonal. Brahmagupta had given a construction 
({\it Br\a hmasphu\d tasiddh\a nta}, XII - 38) of quadrilaterals 
with integer values for the sides and area, starting with a pair of 
Pythagorean triples (the reader is referred to Pranesachar, 2012
\footnote{C.R. Pranesachar, Brahmagupta, mathematician par excellence, 
Resonance, March 2012, 247 - 252.} for an 
exposition on this), and it has been recalled in the {\it L\ndi l\a vat\ndi}
(s\u tras 191-192 : [186-187]). For each quadrilateral constructed using the 
Brahmagupta  construction (which necessarily turn out to be cyclic) 
one gets some new ones (with vertices on the same circle as the original 
one) by replacing the triangle on one side of a 
diagonal by its reflection in the perpendicular bisector of the diagonal. 
This process of obtaining new quadrilaterals also turns up in the 
{\it L\ndi l\a vat\ndi} (though there is no reference to their 
cyclicity along with it). Commentator Ga\d ne\'sa has pointed out 
that the quadrilateral with sides 51, 68, 75 and 40 as in the above 
discussion is one of the quadrilaterals arising in this way, starting
with the Pythagorean triples $(3,4,5)$ and $(8,15,17)$ (cf. Colebrooke, 1993,  
p.~127-128). 

\medskip
\begin{flushleft}
S.G. Dani\\
Department of Mathematics\\
Indian Institute of Technology Bombay \\
Powai, Mumbai 400076\\
India

\smallskip
E-mail :  {\tt sdani@math.iitb.ac.in}
\end{flushleft}

\pagebreak
{

}

\vskip1cm
\begin{flushleft}
S.G. Dani\\
Department of Mathematics\\
Indian Institute of Technology Bombay\\
Powai, Mumbai 400005\\
India

\smallskip
E-mail: {\tt sdani@math.iitb.ac.in}

\end{flushleft}
\end{document}